\documentclass[11pt]{article}
\usepackage{amsthm,amsfonts,amsmath,amssymb}
\newtheorem{thm}{Theorem}

\newtheorem{lem}[thm]{Lemma}

\newtheorem{rem}[thm]{Remark}

\def\di{\bigm|} \def\lg{\langle} \def\rg{\rangle}
\def\nd{\mathrel{\bigm|\kern-.7em/}}

\def\f{\noindent}

\def\mod{\hbox{\rm mod }}

\def\demo{\f{\bf Proof}\hskip10pt}

\def\Diag{\hbox{\rm Diag}}

\def\qed{\hfill $\Box$}

\def\lg{\langle}
\def\rg{\rangle}

\def\B{\mathcal{B}}

\begin{document}
\title{Groups whose Chermak-Delgado lattice is a quasi-antichain\thanks{This work was supported by NSFC (No. 11471198)}}
\author{Lijian An\\
Department of Mathematics,
Shanxi Normal University\\
Linfen, Shanxi 041004, P. R. China\\
 }

\maketitle

\begin{abstract}

A quasiantichain is a lattice consisting of a maximum, a minimum, and the atoms of the lattice. The width of a quasiantichian is the number of atoms.
For a positive integer $w$ ($\ge 3$), a quasiantichain of width $w$ is denoted by $\mathcal{M}_{w}$. In \cite{BHW2}, it is proved that  $\mathcal{M}_{w}$ can be as a Chermak-Delgado lattice of a finite group if and only if $w=1+p^a$ for some positive integer $a$. Let $t$ be the number of abelian atoms in $\mathcal{CD}(G)$. If $t>2$, then, according to \cite{BHW2}, there exists a positive integer $b$ such that $t=p^b+1$. The converse is still an open question. In this paper, we proved that $a=b$ or $a=2b$.

\medskip

\noindent{\bf Keywords}  finite $p$-groups \ \ Chermak-Delgado lattice \ \ quasi-antichain\ \ finitely generated torsion module

\medskip
 \noindent{\it 2000
Mathematics subject classification:} 20D15.
\end{abstract}

\baselineskip=16pt

 Chermak and Delgado \cite{CD} defined a family of functions from the set of subgroups of a finite group into the set of positive integers. They then used these functions to obtain a variety of results, including a proof that every finite group $G$ has a characteristic abelian subgroup $N$ such that $|G : N|\le |G : A|^2$ for any abelian $A\le G$. Chermak-Delgado measures are values of one of these functions.
For any subgroup $H$ of $G$, the Chermak-Delgado measure of $H$ (in $G$) is denoted by $m_G(H)$, and defined as $m_G(H)=|H||C_G(H)|.$
The maximal Chermak-Delgado measure of $G$ is denoted by $m^*(G)$. That is
, $m^*(G)=\max\{ m_G(H)\mid H\le G\}.$
In \cite{I}, Isaacs showed the following theorem.

\begin{thm}{\rm (\cite{I})}\label{th1} Given a finite group $G$, let $\mathcal{CD}(G)=\{ H\mid m_G(H)=m^*(G)\}.$ Then

{\rm (1)} $\mathcal{CD}(G)$ is a lattice of subgroups of $G$.

{\rm (2)} If $H,K\in\mathcal{CD}(G)$, then $\lg H,K\rg=HK$.

{\rm (3)} If $H\in\mathcal{CD}(G)$, then $C_G(H)\in\mathcal{CD}(G)$ and $C_G(C_G(H))=H$.
\end{thm}

By Theorem \ref{th1}, Chermak-Delgado lattice of a finite groups is always a self-dual lattice. It is natural to ask a question: which types of self-dual lattices can be used as Chermak-Delgado lattices of finite groups. Some special cases of the above question are proposed and solved. In \cite{BHW1}, it is proved that, for any integer $n$, a chain of length $n$ can be a Chermak-Delgado lattice of a finite $p$-group.
In \cite{ABQW}, general conclusions are given.

 \begin{thm}{\rm (\cite{ABQW})}\label{th=abqw}
 If $\mathcal{L}$ is a Chermak-Delgado lattice of a finite $p$-group $G$ such that both $G/Z(G)$ and $G'$ are elementary abelian, then are $\mathcal{L}^+$ and $\mathcal{L}^{++}$, where $\mathcal{L}^+$ is a mixed 3-string with center component isomorphic to $\mathcal{L}$ and the remaining components being $m$-diamonds {\rm (}a lattice with subgroups in the configuration of an $m$-dimensional cube{\rm )}, $\mathcal{L}^{++}$ is a mixed 3-string with center component isomorphic to $\mathcal{L}$ and the remaining components being lattice isomorphic to $\mathcal{M}_{p+1}$ {\rm (}a quasiantichain of width $p+1$, see the following definition{\rm )}.
\end{thm}

A quasiantichain is a lattice consisting of a maximum, a minimum, and the atoms of the lattice. The width of a quasiantichian is the number of atoms.
For a positive integer $w$ ($\ge 3$), a quasiantichain of width $w$ is denoted by $\mathcal{M}_{w}$. In \cite{BHW2}, it is proved that  $\mathcal{M}_{w}$ can be as a Chermak-Delgado lattice of a finite group if and only if $w=1+p^a$ for some positive integer $a$. According to \cite{BHW2}, if $\mathcal{M}_{w}$ is a Chermak-Delgado lattice of a finite group $G$, then $G$ is nilpotent of class 2. Moreover, we may choose $G$ to be a $p$-group of class 2 without loss of generality. Let $M_1,M_2,\dots,M_w$ be all atoms of $\mathcal{CD}(G)$. Then $G$ has the following properties:

{\rm (P1)} Both $G/Z(G)$ and $G'$ are elementary abelian;

{\rm (P2)} $M_i/Z(G)\cong M_j/Z(G)$ and $G/Z(G)=M_i/Z(G)\times M_j/Z(G)$ for $i\ne j$;

{\rm (P3)} Let $M_1=\lg x_1,x_2,\dots,x_n\rg Z(G)$ and $M_2=\lg y_1,y_2,\dots,y_n\rg Z(G)$ such that $M_1/Z(G)$ and $M_2/Z(G)$ are elementary abelian groups of order $p^n$. Then, for $k\ge 3$, $M_k=\lg x_1y_1',x_2y_2',\dots,x_ny_n'\rg Z(G)$ where $M_2=\lg y_1',y_2',\dots,y_n'\rg Z(G)$. Moreover, there is an invertible matrix $C_k=(c^{(k)}_{ij})_{n\times n}$ over $F_p$ such that
$$y_j'\equiv \prod_{i=1}^{n} y_i^{c^{(k)}_{ij}} \ (\mod Z(G)), j=1,2,\dots,n.$$
In this paper, $C_k$ is called the characteristic matrix of $M_k$ relative to $(x_1,x_2,\dots,x_n)$ and $(y_1,y_2,\dots,y_n)$.

Let $G$ be a $p$-group with $\mathcal{CD}(G)$ a quasi-antichain of width $w\ge 3$. Let $t$ be the number of abelian atoms in $\mathcal{CD}(G)$. If $t>2$, then, according to \cite{BHW2}, there exists a positive integer $b$ such that $t=p^b+1$. The converse is still an open question. At the end of \cite{BHW2}, such questions are proposed definitely:

{\rm (Q1)} Which values of $t$ are possible in quasi-antichain Chermak-Delgado lattices of width $w=p^a+1$ when $a>1$?

{\rm (Q2)} Are there examples of groups $G$ with $G\in \mathcal{CD}(G)$ and $\mathcal{CD}(G)$ a quasi-antichain where $t=0$ and $p\equiv 1$ modulo $4$?

This paper answer above questions completely. It is amazing that there are only two possible relations: $a=b$ or $a=2b$. We also prove that $a\mid n$, which is another  application about the decomposition of a finitely generated torsion module over a principal ideal domain. Examples related to (Q1) and (Q2) are given.

Throughout the article, let $p$ be a prime and $G$ be a finite $p$-group with properties (P1), (P2) and (P3).

\begin{thm}\label{thm=main}
Suppose that there are totally $t>2$ abelian atoms in $\mathcal{{CD}}(G)$, which is a quasi-antichain of width $w$. Then there are positive integers $a$ and $b$ such that $w=p^a+1$ and $t=p^b+1$, where $a=b$ or $a=2b$. If $|G/Z(G)|=p^{2n}$, then $a\mid n$.
\end{thm}
{\bf Proof}\quad Without loss of generality, we may assume that $M_1$, $M_2$ and $M_3$ are abelian atoms.
 For convenience, the operation of the group $G$ is replaced to be addition. According to (P3), $M_3=\lg x_1+y_1',x_2+y_2',\dots,x_n+y_n'\rg Z(G)$ such that
$$y_j'\equiv \prod_{i=1}^{n} {c^{(3)}_{ij}y_i} \ (\mod Z(G)), j=1,2,\dots,n.$$
 where $C_3=(c^{(3)}_{ij})_{n\times n}$ is the characteristic matrix of $M_k$ relative to $(x_1,x_2,\dots,x_n)$ and $(y_1,y_2,\dots,y_n)$. Replacing $y_1,y_2, \dots, y_n$ with $y'_1,y'_2,\dots,y'_n$ respectively, we may assume that $C_3=I_n$.

%
%In this article, the operation of the group $G$ is replaced to be addition. According to (P3), for $k\ge 3$, we may let $M_k=\lg w_1,w_2,\dots,w_n\rg Z(G)$ such that
%$$w_i=x_i+\prod_{j=1}^{n} {c_{ij}y_j}, i=1,2,\dots,n.$$
%where $C_k=(c_{ij})_{n\times n}$ is a matrix over $F_p$ (the
%finite field containing $p$ elements).

 Let $z_{ij}=[x_i,y_j]$ and $Z=(z_{ij})_{n\times n}$. Then $G'=\lg z_{ij}\mid i,j=1,2,\dots,n\rg$.
 We have the following formalized calculation:
\begin{eqnarray*}
  &[((x_1,x_2,\dots,x_n)+(y_1,y_2,\dots,y_n)C_k)^T,(x_1,x_2,\dots,x_n)+(y_1,y_2,\dots,y_n)C_{k'}]  \\
  =&[(x_1,x_2,\dots,x_n)^T+C_k^T(y_1,y_2,\dots,y_n)^T,(x_1,x_2,\dots,x_n)+(y_1,y_2,\dots,y_n)C_{k'}]\\
  =&[(x_1,x_2,\dots,x_n)^T,(y_1,y_2,\dots,y_n)C_{k'}]+[C_k^T(y_1,y_2,\dots,y_n)^T,(x_1,x_2,\dots,x_n)]\\
  =& ZC_{k'}+C_k^T(-Z^T)=ZC_{k'}-(ZC_k)^T.
\end{eqnarray*}
Hence $C_G(M_k)=M_{k'}$ if and only if $(ZC_k)^T=ZC_{k'}$. In particular, $M_k$ is abelian if and only if $ZC_{k}$ is symmetric. Since $M_3$ is abelian, $Z^T=Z$.

Let $\mathbb{V}=\{ C\in \mathcal{M}_n(F_p)\mid \exists  \ C'\in \mathcal{M}_n(F_p)\  {\rm s.t.} \ C^TZ=ZC'\}$. It is straightforward to prove that $\mathbb{V}$ is a linear space over $F_p$. Hence $|\mathbb{V}|=p^a$ for some positive integer $a$. In the following, we prove that $\mathbb{V}$ is also a field.

Assume that $O_n\ne C\in\mathbb{V}$ and $C^TZ=ZC'$. Let $$(s_1,s_2,\dots,s_n)=(x_1,x_2,\dots,x_n)+(y_1,y_2,\dots,y_n)C$$ and $$(t_1,t_2,\dots,t_n)=(x_1,x_2,\dots,x_n)+(y_1,y_2,\dots,y_n)C'.$$ Then, by the above formalized calculation, $[(s_1,s_2,\dots,s_n)^T,(t_1,t_2,\dots,t_n)]=O_n$. Let $N_1=\lg s_1,s_2,\dots,s_n\rg Z(G)$
 and $N_2=\lg t_1,t_2,\dots,t_n\rg Z(G)$. Then $N_2\le C_G(N_1)$. Hence $|N_1||C_G(N_1)|\ge |N_1||N_2|=|M_1|^2=m^*(G)$. It implies that $N_1\in \mathcal{CD}(G)$, and, there is some $k\ge 3$ such that $N_1=M_k$.
 Hence $\mathbb{V}\setminus \{ O_n\}$ is the set of $C_k$ where $3 \le k\le w$. Since $C_k$ is invertible, $\mathbb{V}$ is a division algebra. By Wedderburn's little theorem, $\mathbb{V}$ is also a finite field.
 Now we know that $w-2=p^a-1$. hence $w=p^a+1$.  Let $A$ be a primitive element of $\mathbb{V}$. Then $\mathbb{V}=\{ O_n,A,A^2,\dots,A^{p^a-1}=I_n\}$.

Next we study $A$ applying decomposition of a finitely generated torsion module over a principal ideal domain.
$M_2/Z(G)$ can be regarded as a $n$-dimensional vector space $\mathbb{V}_2$ over $F_p$, in which $(y_1Z(G),y_2Z(G),\dots,y_nZ(G))$ is a base for $\mathbb{V}_2$. At this time, $A=(a_{ij})$ is a matrix of a single transformation $\mathcal{A}$ relative to the given base, where
$$\mathcal{A}( y_jZ(G))=(\prod_{i=1}^{n} {a_{ij}y_i) Z(G)} , j=1,2,\dots,n.$$
We can make $\mathbb{V}_2$ an $F_p[\lambda]$-module by defining the action of any polynomial $g(\lambda)=b_0+b_1\lambda+\dots+b_m\lambda^m$ on any vector $y\in \mathbb{V}_2$ as
$$g(\lambda)y=b_0+b_1(\mathcal{A}y)+\dots+b_m(\mathcal{A}^my).$$
Let $m_A(\lambda)$ be the minimal polynomial of $A$. Then $\mathbb{V}_2$ is a torsion $F_p[\lambda]$-module of exponent $m_A(\lambda)$. Obviously, $m_A(x)$ is irreducible with degree $a$.
By fundamental structure theorem for finitely generated torsion modules over a principle ideal domain, $\mathbb{V}_2$ is a direct sum of cyclic module:
$$\mathbb{V}_2=F_p[\lambda] v_1\oplus F_p[\lambda] v_2\oplus\dots\oplus F_p[\lambda]v_r$$
It is easy to see that $ann v_i=m_A(\lambda)$ where $i=1,2,\dots,r$. Let
$m_A(\lambda)=k_0+k_1\lambda+\dots+k_{a-1}\lambda^{a-1}-k_a\lambda^{a}$. Then
$n=ra$ and
$$v_1,\mathcal{A}v_1,\dots,\mathcal{A}^{a-1}v_1,v_2,\mathcal{A}v_2,\dots,\mathcal{A}^{a-1}v_2,\dots,v_r,\mathcal{A}v_r,\dots,\mathcal{A}^{a-1}v_r$$
is a set of generators for $M_2$.
Assume that $(v_1, v_2,\dots,v_n)=(y_1,y_2,\dots,y_n)S$. Let $(u_1, u_2 ,\dots u_n)=(x_1,x_2,\dots,x_n)S$
 where $S=(s_{ij})$ is an invertible matrix. Then the characteristic matrix of $M_k$ relative to $(u_1,u_2,\dots,u_n)$ and $(v_1, v_2,\dots,v_n)$ is $S^{-1}C_kS$.
 Matrices related in this way are said to be similar.

Without loss of generality, we may assume that $A$ is the characteristic matrix relative to $(v_1, \mathcal{A}v_1,\dots,\mathcal{A}^{a-1}v_r)$. It is clear that $A$ has the form
  $$A=\left(
                                         \begin{array}{cccc}
                                           A_1 & O_a & \dots & O_a \\
                                           O_a & A_2 & \dots & O_a    \\
                                           \vdots &\vdots&\ddots&\vdots \\
                                           O_a & 0 & \dots & A_r\\
                                         \end{array}
                                       \right),$$
where $A_i$ is the companion matrix of $m_A(\lambda)$. That is,
$$A_i=B=\left(
                                         \begin{array}{cccc}
                                           0 & 0  & \dots & k_0 \\
                                           1 & 0  & \dots & k_1 \\
                                          % 0 & 1 & 0 & \dots & k_2 \\
                                           \vdots &  \ddots &\ddots&\vdots \\
                                            0 & \dots& 1 & k_{a-1}\\
                                         \end{array}
                                       \right),$$
and $A=\Diag(B,B,\dots, B)$.

If $A^TZ=ZA'$, then $A'TZ=ZA$ since $Z^T=Z$. It follows that $A'\in\mathbb{V}$. Hence there are $1\le k\le p^a-1$ such that $A^TZ=ZA^k$. Let
 $$Z=\left(
                                         \begin{array}{cccc}
                                           Z_{11} & Z_{12} & \dots & Z_{1r} \\
                                            Z_{21} & Z_{22} & \dots & Z_{2r} \\
                                           \vdots &\vdots&\ddots&\vdots \\
                                            Z_{r1} & Z_{r2} & \dots & Z_{rr}\\
                                         \end{array}
                                       \right),$$
 where $Z_{ij}$ are $a\times a$ matrices. Since $Z^T=Z$, $Z_{ji}=Z_{ij}^T$ for $1\le i\le j\le r$. It follows from $A^TZ=ZA^k$ that
 \begin{eqnarray}
 \label{eq:11}  B^TZ_{ij} &=& Z_{ij}B^k \\
   B^TZ_{ji} &=& Z_{ji}B^k \label{eq:12}
 \end{eqnarray}
 Let $Z_{ij}=\left(\begin{array}{c}
                                       Z_{ij}^{(1)} \\
                                       Z_{ij}^{(2)} \\
                                       \vdots\\
                                       Z_{ij}^{(a)} \\
                                     \end{array}
                                   \right)$. Then
                                   \begin{equation}\label{eq:13}
                                   B^TZ_{ij}=\left(\begin{array}{c}
                                       Z_{ij}^{(2)} \\
                                       Z_{ij}^{(3)} \\
                                       \vdots\\
                                       k_0 Z_{ij}^{(1)}+k_1Z_{ij}^{(2)}+\dots+k_{a-1}Z_{ij}^{(a)} \\
                                     \end{array}
                                   \right)
 \end{equation}
 By (\ref{eq:11}) and (\ref{eq:13}),
  \begin{equation}\label{eq:14}
 Z_{ij}^{(2)}=Z_{ij}^{(1)}B^k, Z_{ij}^{(3)}=Z_{ij}^{(2)}B^k, \dots,Z_{ij}^{(a)}=Z_{ij}^{(a-1)}B^k
 \end{equation}
  \begin{equation}\label{eq:15}
 k_0Z_{ij}^{(1)}+k_1Z_{ij}^{(2)}+\dots+k_{a-1}Z_{ij}^{(a)}=Z_{ij}^{(a)}B^k
 \end{equation}
Hence
  \begin{equation}\label{eq:16}
 Z_{ij}^{(1)}(k_0+k_1B^k+\dots+k_{a-1}B^{(a-1)k}-B^{ak})=0
 \end{equation}
 That is, $Z_{ij}^{(1)}m_A(B^k)=0$. We claim that $m_A(B^k)=0$. Otherwise, $m_A(A^k)\ne 0$. Since $m_A(A^k)\in \mathbb{V}$ and $\mathbb{V}$ is a field, $m_A(A^k)$ is invertible. Hence $m_A(B^k)$ is also invertible.
Notice that if $Z_{ij}^{(1)}=0$ then $Z_{ij}=O_a$ by (\ref{eq:14}) and (\ref{eq:15}).
Hence we may choose $Z_{ij}$ such that $Z_{ij}^{(1)}\ne 0$. In this case, $Z_{ij}^{(1)}m_A(B^k)\ne 0$, a contradiction. Thus $m_A(B^k)=0$ and $m_A(A^k)=0$. Since $A^p,A^{p^2},\dots,A^{p^a}=A$ are all zero points of $m_A(\lambda)$, there exists a $1\le e\le a$ such that $k=p^e$.

 It is easy to see that $(A^m)^TZ=ZA^{mp^e}$ if and only if $(p^{a}-1)\mid m(p^e-1)$. Let $\mathbb{W}=\{ C\in \mathbb{V}\mid C^TZ=ZC \}.$ Then $\mathbb{W}\setminus\{ O_n\}=\{ A^m\in \mathbb{V}\mid (p^{a}-1)\mid m(p^e-1) \}.$ Hence $\mathbb{W}\setminus\{ O_n\}$ is a cyclic group generated by $A^{\frac{p^a-1}{(p^a-1,p^e-1)}}$ of order $(p^a-1,p^e-1)=p^{(a,e)}-1$. Let $b=(a,e)$. Notice that $|\mathbb{W}|=t-1$. It is easy to see that $t=p^b+1$.
By (\ref{eq:11}),
\begin{equation}\label{eq:17}
   (B^{p^e})^TZ_{ij} =Z_{ij}B^{p^{2e}}
\end{equation}
 By (\ref{eq:12}),
 \begin{equation}\label{eq:18}
   (B^{p^e})^TZ_{ij} =Z_{ij}B
 \end{equation}
 Hence $(p^a-1)\mid (p^{2e}-1)$. Thus $a\mid 2e$. It follows that $a\mid 2b$. Hence $a=e=b$ or $a=2e=2b$.

\qed

\begin{rem}\label{rem} Let $F$ be a field containing $p^a$ elements and $F^*$ be the
multiply group of $F$. Then $F^*$ is cyclic with order $p^a-1$. Let
$F^*=\lg b\rg$ and $\B: f\mapsto bf$ be linear transformation over
$F$, where $F$ is regarded as a linear space over field $F_p$. Of
course, the order of $\B$ is $p^a-1$. Let $p(x)$ be the minimal
polynomial of $\B$. It follows from Cayley-Hamilton theorem that
$deg(p(x))=r\le a$. Let $W=\{f(\B)\di f(x)\in F_p[x]\}=\{f_1(\B)\mid
f_1(x)\in F_p[x], deg(f_1)<r\}$. Then $dim W=deg(p(x))=r$ and hence
$|W|=p^r$. On the other hand,
$\{1,\B,\B^2,\dots,\B^{p^a-1}\}\subseteq W$ and hence $|W|\ge p^a$.
So we get $r=a$. Let the minimal polynomial of $\B$ is
$$p(x)=x^a-k_{a-1}x^{a-1}-\dots-k_1x-k_0.$$ Then a matrix of $\B$ is
Frobenius form $$B=\left(
                                         \begin{array}{cccc}
                                           0 & 0  & \dots & k_0 \\
                                           1 & 0  & \dots & k_1 \\
                                          % 0 & 1 & 0 & \dots & k_2 \\
                                           \vdots &  \ddots &\ddots&\vdots \\
                                            0 & \dots& 1 & k_{a-1}\\
                                         \end{array}
                                       \right).$$
\end{rem}

\begin{lem}\label{lem=ss} Let $B$ be the matrix introduced in Remark $\ref{rem}$. $Z=(z_{ij})$ is a $a\times a$ matrix. If $B^TZ=ZB$, Then $Z$ is symmetric.
\end{lem}
\demo By calculation, the $(u,v)^{th}$ entry of $B^TZ$ is $z_{u+1,v}$, while the $(u,v)^{th}$ entry of $ZB$ is $z_{u,v+1}$, where $1\le u,v\le a-1$. Hence $z_{u+1,v}=z_{u,v+1}$ where $1\le u,v\le a-1$. Then, for $u\le v$, we have $z_{u,v+1}=z_{u+1,v}=z_{u+2,v-1}=\dots=z_{v,u+1}=z_{v+1,u}$. Hence $Z$ is symmetric.

\begin{thm}\label{thm=ex1}
 Suppose that ???
\end{thm}

\demo
Let $G$ be generated by $\{x_1,x_2,\dots,x_{n},y_1,y_2,\dots,y_{n}\}$ with defining relationships $x_i^p=y_i^p=1$ and $[x_i,x_j]=[y_i,y_j]=1$ for all $i,j$ such that $1\le i,j\le n$, $[x_{(u-1)a+i},y_{(v-1)a+j}]=z_{(u-1)a+i,(v-1)a+j}$ for $1\le u,v \le r$ and $1\le i,j\le a$,
% $\lg z_{1+ia,j}\mid 0\le i\le r-1, 1+ia\le j\le n\rg\le Z(G)$ is elementary abelian with order $p^{\frac{r+1}{2}n}$.
$z_{(u-1)a+1,(v-1)a+j}\in Z(G)$ for $1\le u\le v \le r$ and $1\le j\le a$.
 For convenience, we use addition operation to replace mutiplication operation of $G$. We also use the following notations (where $1\le u,v\le r$ and $1\le i,j,k\le a$):
$$Z_{uv}^{(k)}:=(z_{(u-1)a+k,(v-1)a+1},z_{(u-1)a+k,(v-1)a+2},\dots,z_{(u-1)a+k,(v-1)a+a}),$$
$$Z_{uv}=(z_{(u-1)a+i,(v-1)a+j})=\left(\begin{array}{c}
                                       Z_{uv}^{(1)} \\
                                      Z_{uv}^{(2)}\\
                                       \vdots\\
                                      Z_{uv}^{(3)}  \\
                                     \end{array}
                                   \right)\ \mbox{and}\ Z=\left(
                                         \begin{array}{cccc}
                                           Z_{11} & Z_{12} & \dots & Z_{1r} \\
                                            Z_{21} & Z_{22} & \dots & Z_{2r} \\
                                           \vdots &\vdots&\ddots&\vdots \\
                                            Z_{r1} & Z_{r2} & \dots & Z_{rr}\\
                                         \end{array}
                                       \right).$$
Using above notations, we continue to give defining relationships $Z_{uv}^{(k)}=Z_{uv}^{(1)}B^{k-1}$ for $1\le u\le v\le r$ and $2\le k\le a$, $Z_{uv}=Z_{vu}^T$ for $1\le v<u\le r$. It is easy to see that
 $G'=Z(G)=\lg z_{(u-1)a+1,(v-1)a+j}\mid 1\le u\le v \le r, 1\le j\le a\rg$ is elementary abelian of order $p^{\frac{r+1}{2}n}$. Hence $|G|=p^{\frac{r+5}{2}n}$.

Since $Z_{uv}^{(k)}=Z_{uv}^{(1)}B^{k-1}$ for $1\le u\le v\le r$ and $2\le k\le a$, $B^TZ_{uv}=Z_{uv}B$ for $1\le u\le v\le r$. By Lemma \ref{lem=ss}, $Z_{uv}$ is symmetric for $1\le u\le v\le r$. Hence $Z_{uv}=Z_{vu}^T=Z_{vu}$ for $1\le v<u\le r$. Moreover $Z_{uv}^T=Z_{vu}^T=Z_{uv}=Z_{vu}$ for all $1\le u,v\le r$ and $Z^T=Z$.
 Let $X=\lg x_1,x_2,\dots,x_{n}\rg Z(G)$ and $Y=\lg y_1,y_2,\dots,y_{n}\rg Z(G)$.

\medskip

Assertion 1:  $C_G(x)=X$ for all $x\in X\setminus Z(G)$.

Let $x=\prod_{i=1}^{n}c_{i}x_i+z$ where $z\in Z(G)$. Write $C_k=(c_{(k-1)a+1},c_{(k-1)a+2},\dots,c_{ka})$ for $1\le k\le a$. Since $x\not\in Z(G)$, there exists a $k_0$ such that $C_{k_0}\ne (0,0,\dots,0)$.
Exchanging $x_{(k_0-1)a+1},y_{(k_0-1)a+1},\dots,x_{k_0a},y_{k_0a}$ and $x_{(r-1)a+1},y_{(r-1)a+1},\dots,x_{ra},y_{ra}$, we have $C_r\ne (0,0,\dots,0)$.
Let $H_i=\lg [x,y_{(i-1)a+j}]\mid 1\le j\le a\rg$ for $1\le i\le r$.
We will prove that $H_1+H_2+\dots+H_v$ is of order $p^{va}$. Use induction, we may assume that $H_1+H_2+\dots+H_{v-1}$ is of order $p^{(v-1)a}$. By calculation,
\begin{eqnarray*}
 & &([x,y_{(v-1)a+1}],\dots,[x,y_{va}])\\
 &=& \sum_{k=1}^n c_k([x_k,y_{(v-1)a+1}],\dots,[x_k,y_{va}]) \\
&=& \sum_{u=1}^r \sum_{i=1}^a c_{(u-1)a+i}([x_{(u-1)a+i},y_{(v-1)a+1}],\dots,[x_{(u-1)a+i},y_{va}])\\
&=& \sum_{u=1}^r \sum_{i=1}^a c_{(u-1)a+i}Z_{uv}^{(i)}\\
&=& \sum_{u=1}^r \sum_{i=1}^a c_{(u-1)a+i}Z_{uv}^{(1)}B^{i-1}\\
&=& \sum_{u=1}^r Z_{uv}^{(1)}(\sum_{i=1}^a  c_{(u-1)a+i}B^{i-1})\\
&=& \sum_{u=1}^{r-1} Z_{uv}^{(1)}(\sum_{i=1}^a  c_{(u-1)a+i}B^{i-1})+Z_{rv}^{(1)}(\sum_{i=1}^a  c_{(r-1)a+i}B^{i-1}) \\
\end{eqnarray*}
Since $C_r\ne (0,0,\dots,0)$, $(\sum_{i=1}^a  c_{(r-1)a+i}B^{i-1})\ne O_a$. By Remark \ref{rem}, $(\sum_{i=1}^a  c_{(r-1)a+i}B^{i-1})$ is invertible. Hence $Z_{rv}^{(1)}(\sum_{i=1}^a c_{(r-1)a+i}B^{i-1})$ is of rank $a$. It follows that $(H_1+H_2+\dots+H_{v-1})\cap H_v=0$. Thus $H_1+H_2+\dots+H_v$ is of order $p^{va}$.

By above discussion, $H_1+H_2+\dots+H_v$ is of order $p^{n}$.
Hence $|[x,G]|=p^{n}$. It follows that $|C_G(x)|=|G|/p^{n}=p^{\frac{r+3}{2}n}$.
Since $|X|=p^{\frac{r+3}{2}n}$ and $X\le C_G(x)$, $C_G(x)=X$.

Similarly,  we have $C_G(y)=Y$ for all $y\in Y\setminus Z(G)$.
Then $C_G(X)=X$ and $C_G(Y)=Y$, yielding $m_G(G)=m_G(X)=m_G(Y)=p^{(r+3)n}$.

\medskip

Assertion 2: $m^*(G)=p^{(r+3)n}$.

Otherwise, by the dual-property of $\mathcal{CD}$-lattice, there exists $H\in
\mathcal{CD}(G)$ such that $H<G$ and $|H|>p^{\frac{(r+3)}{2}n}$. Since $|H\cap
X|=\frac{|H||X|}{|HX|}>p^{\frac{(r+1)}{2}n}=|Z(G)|$, there exists $x\in H\cap X\setminus
Z(G)$. Hence $C_G(H)\le C_G(x)=X$. Similarly, we have $C_G(H)\le Y$. Hence $C_G(H)=Z(G)$ and $m_G(H)<m_G(G)$, a contradiction.

Above discussion also gives that $\mathcal{CD}(G)$ is a quasi-antichain, in which every atom is of order $p^{\frac{(r+3)}{2}n}$.

\medskip

Assertion 3:  $w=t=1+p^a$.

Now we have $G,Z(G),X,Y\in \mathcal{CD}(G)$. Let $M$ be an atom different from $X$ and $Y$. Then, by the same reason given in Theorem \ref{thm=main}, we may let $M=\lg w_1,w_2,\dots,w_{n}\rg$ and $N=C_G(M)=\lg v_1,v_2,\dots,v_{n}\rg$ where
$$(w_1,w_2,\dots,w_{n})=(x_1,x_2,\dots,x_{n})+(y_1,y_2,\dots,y_{n})C,$$
$$(v_1,v_2,\dots,v_{n})=(x_1,x_2,\dots,x_{3n})+(y_1,y_2,\dots,y_{n})D.$$
Since $[M,N]=0$, $C^TZ=ZD$. Let
$$C=\left(
                                         \begin{array}{cccc}
                                           C_{11} & C_{12} & \dots & C_{1r} \\
                                            C_{21} & C_{22} & \dots & C_{2r} \\
                                           \vdots &\vdots&\ddots&\vdots \\
                                            C_{r1} & C_{r2} & \dots & C_{rr}\\
                                         \end{array}
                                       \right)\ \mbox{and}\ D=\left(
                                         \begin{array}{cccc}
                                           D_{11} & D_{12} & \dots & D_{1r} \\
                                            D_{21} & D_{22} & \dots & D_{2r} \\
                                           \vdots &\vdots&\ddots&\vdots \\
                                            D_{r1} & D_{r2} & \dots & D_{rr}\\
                                         \end{array}
                                       \right).$$
Then
\begin{equation}\label{eq:19}
\sum_{k=1}^r C_{ki}^TZ_{kj}=\sum_{k=1}^r Z_{ik} D_{kj} \ \mbox{for} \ 1\le i,j\le r.
\end{equation}
Taking $i\ne j$ in Equation (\ref{eq:19}) and compiling two sides, we get $C_{ik}=O_a$ for $i\ne k$, $D_{kj}=O_a$ for $k\ne j$, and
\begin{equation}\label{eq:110}
C_{ii}^TZ_{ij}=Z_{ij} D_{jj} \ \mbox{for} \ 1\le i,j\le r \ \mbox{and} \ i\ne j.
\end{equation}
Taking $i=j$ in Equation (\ref{eq:19}) and compiling two sides, we get
\begin{equation}\label{eq:111}
C_{ii}^TZ_{ii}=Z_{ii} D_{ii} \ \mbox{for} \ 1\le i\le r.
\end{equation}
Notice that $Z_{ij}$ and $Z_{11}$ have no essential difference. Thus $C_{ii}=C_{jj}$ and $D_{ii}=D_{jj}$ for $i\ne j$.

Let $\mathbb{V}=\{ U\in \mathcal{M}_a(F_p)\mid \exists  \ V\in \mathcal{M}_a(F_p)\  {\rm s.t.} \ U^TZ_{11}=Z_{11}V\}$.
Similar to the proof in Theorem \ref{thm=main}, $\mathbb{V}$ is a division algebra. Hence $\mathbb{V}$ has at most $p^a$ elements (otherwise, there exists a non-zero matrix in which the elements of the first row are all zero, a contradiction). Notice that $B^TZ_{11}=Z_{11}B$. Hence $B^k\in \mathbb{V}$ where $1\le k\le p^a-1$. it follows that $\mathbb{V}=\{ O_n,B,B^2,\dots,B^{p^a-1}=I_n\}$.

Therefore $C=D=\Diag(B^k,B^k,\dots,B^k)$ for some $1\le k\le p^a-1$. Hence $w=t=p^a+1$.
\qed

\begin{thm}\label{thm=ex2}
 Suppose that $a,b$ and $r$ are positive integers such that $a=2b$ and $r\ge 3$. Let $n=ar$. Then there exists a group $G$ such that $|G/Z(G)|=p^{2n}$, and $\mathcal{CD}(G)$ is a quasi-antichain of width $w=p^a+1$, in which the number of abelian atoms is $t=p^b+1$.
\end{thm}

\demo
Let $G$ be generated by $\{x_1,x_2,\dots,x_{n},y_1,y_2,\dots,y_{n}\}$ with defining relationships similar to that in Theorem \ref{thm=ex1}. The followings are different defining relationships from that in Theorem \ref{thm=ex1}:
$[x_{(u-1)a+i},y_{(u-1)a+j}]=1$ for $1\le u \le r$ and $1\le i,j\le a$, $Z_{uv}^{(k)}=Z_{uv}^{(1)}B^{(k-1)p^b}$ for $1\le u<v\le r$ and $2\le k\le a$.

In this case,
 $G'=Z(G)=\lg z_{(u-1)a+1,(v-1)a+j}\mid 1\le u<v \le r, 1\le j\le a\rg$ is elementary abelian of order $p^{\frac{r-1}{2}n}$. Hence $|G|=p^{\frac{r+3}{2}n}$.

Since $Z_{uv}^{(k)}=Z_{uv}^{(1)}B^{(k-1)p^b}$ for $1\le u<v\le r$ and $2\le k\le a$, $B^TZ_{uv}=Z_{uv}B^{p^b}$ for $1\le u<v\le r$. In this case, $Z_{uv}$ is not symmetric and $Z_{uv}\ne Z_{vu}$ for all $1\le v<u\le r$. $Z^T=Z$ is still hold.
 Let $X=\lg x_1,x_2,\dots,x_{n}\rg Z(G)$ and $Y=\lg y_1,y_2,\dots,y_{n}\rg Z(G)$.

Similar to Assertion 1 in the proof of Theorem \ref{thm=ex1}, we have $|[x,Y]|\ge p^{(r-1)a}$ for all $x\in X\setminus Z(G)$. It follows that $|C_Y(x)|\le |Y|/p^{(r-1)a}=p^a|Z(G)|$.

Similarly, $|C_X(y)|\le p^{a}|Z(G)|$ for all $y\in Y\setminus Z(G)$.
It is easy to check that $C_G(X)=X$ and $C_G(Y)=Y$, yielding $m_G(G)=m_G(X)=m_G(Y)=p^{(r+1)n}$.

\medskip

Assertion (1): $m^*(G)=p^{(r+1)n}$.

Otherwise, by the dual-property of $\mathcal{CD}$-lattice, there exists $H\in
\mathcal{CD}(G)$ such that $H<G$ and $|H|>p^{\frac{(r+1)}{2}n}$. Since $|H\cap
X|=\frac{|H||X|}{|HX|}>p^{\frac{(r-1)}{2}n}=|Z(G)|$, there exists $x\in H\cap X\setminus
Z(G)$. Hence $C_G(H)\le C_G(x)=XC_Y(x)$. Similarly, there exists $y\in H\cap Y\setminus
Z(G)$. Hence $C_G(H)\le C_G(y)=YC_X(y)$. It follows that $C_G(H)\le C_X(y)C_Y(x)$. Hence $|C_G(H)|\le p^{2a}|Z(G)|$. Obviously $C_G(H)>Z(G)$. Hence there exists $1\ne x'y'\in C_G(H)\setminus Z(G)$, where $x'\in C_X(y), y'\in C_Y(x)$. Whatever, we have $|H|\le |C_G(x'y')|\le |G:[G,x'y']|\le p^{\frac{r+1}{2}n+a}$. Hence $m_G(H)=|H||C_G(H)|\le p^{rn+3a}\le p^{(r+1)n}$,
a contradiction.

\medskip

Assertion (2): $\mathcal{CD}(G)$ is a quasi-antichain.

Otherwise, by the dual-property of $\mathcal{CD}$-lattice, there exists $H\in
\mathcal{CD}(G)$ such that $H<G$ and $|H|>p^{\frac{(r+1)}{2}n}$. By Assertion (1), $n=3a$ and there exists $x\in H\cap X\setminus Z(G)$ and $y\in H\cap Y\setminus Z(G)$ such that $C_G(H)=C_X(y)C_Y(x)$, where $|C_X(y)|=|C_Y(x)|=p^a|Z(G)|$. Take $x'\in C_X(y)\setminus Z(G)$ and $y'\in C_Y(x)\setminus Z(G)$. Then $H\le C_X(y')C_X(x')\le p^{2a}|Z(G)|<p^{\frac{(r+1)}{2}n}$, a contradiction.

\medskip

Assertion (3):  $w=p^a+1$ and $t=p^b+1$.

Let $M$ be an atom different from $X$ and $Y$. Then, by the same reason given in Theorem \ref{thm=main}, we may let $M=\lg w_1,w_2,\dots,w_{n}\rg$ and $N=C_G(M)=\lg v_1,v_2,\dots,v_{3n}\rg$ where
$$(w_1,w_2,\dots,w_{n})=(x_1,x_2,\dots,x_{n})+(y_1,y_2,\dots,y_{n})C,$$
$$(v_1,v_2,\dots,v_{3})=(x_1,x_2,\dots,x_{3n})+(y_1,y_2,\dots,y_{3})D.$$
We also have $C^TZ=ZD$. Let
$$C=\left(
                                         \begin{array}{cccc}
                                           C_{11} & C_{12} & \dots & C_{1r} \\
                                            C_{21} & C_{22} & \dots & C_{2r} \\
                                           \vdots &\vdots&\ddots&\vdots \\
                                            C_{r1} & C_{r2} & \dots & C_{rr}\\
                                         \end{array}
                                       \right)\ \mbox{and}\ D=\left(
                                         \begin{array}{cccc}
                                           D_{11} & D_{12} & \dots & D_{1r} \\
                                            D_{21} & D_{22} & \dots & D_{2r} \\
                                           \vdots &\vdots&\ddots&\vdots \\
                                            D_{r1} & D_{r2} & \dots & D_{rr}\\
                                         \end{array}
                                       \right).$$
Then
\begin{equation}\label{eq:112}
\sum_{k=1}^r C_{ki}^TZ_{kj}=\sum_{k=1}^r Z_{ik} D_{kj} \ \mbox{for} \ 1\le i,j\le r.
\end{equation}
Taking $i\ne j$ in Equation (\ref{eq:112}) and compiling two sides, we get $C_{ik}=O_a$ for $i\ne k$, $D_{kj}=O_a$ for $k\ne j$, and
\begin{equation}\label{eq:110}
C_{ii}^TZ_{ij}=Z_{ij} D_{jj} \ \mbox{for} \ 1\le i,j\le r \ \mbox{and} \ i\ne j.
\end{equation}
Notice that $Z_{ij}$, $Z_{ik}$ and $Z_{ki}$ have no essential difference for $k\ne i,j$. Thus $C_{ii}=C_{kk}$ and $D_{jj}=D_{kk}$ for $k\ne i,j$. Since $r\ge 3$, $C_{11}=C_{22}=\dots=C_{rr}$ and $D_{11}=D_{22}=\dots=D_{rr}$.

Let $\mathbb{V}=\{ U\in \mathcal{M}_a(F_p)\mid \exists  \ V\in \mathcal{M}_a(F_p)\  {\rm s.t.} \ U^TZ_{12}=Z_{12}V\}$.
Similar to the proof in Theorem \ref{thm=main}, $\mathbb{V}$ is a division algebra. Hence $\mathbb{V}$ has at most $p^a$ elements (otherwise, there exists a non-zero matrix in which the elements of the first row are all zero, a contradiction). Notice that $B^TZ_{12}=Z_{12}B^{p^b}$. Hence $B^k\in \mathbb{V}$ where $1\le k\le p^a-1$. It follows that $\mathbb{V}=\{ O_n,B,B^2,\dots,B^{p^a-1}=I_n\}$.

Therefore $C=\Diag(B^k,B^k,\dots,B^k)$ and $D=\Diag(B^{kp^b},B^{kp^b},\dots,B^{kp^b})$ for some $1\le k\le p^a-1$. Hence $w=p^a+1$. It is easy to see that $M=N$ if and only if $C=D$ if and only if $(p^b+1)\mid k$. Hence the number of abelian atoms $t=p^b+1$.
\qed

\begin{thm}\label{thm=ex3}
 Suppose that $a$ and $r$ are positive integers such that $r\ge 3$. Let $n=ar$. Then there exists a group $G$ such that $|G/Z(G)|=p^{2n}$, and $\mathcal{CD}(G)$ is a quasi-antichain of width $w=p^a+1$, in which the number of abelian atoms is $1$ or $2$ for $p=2$ or $p\ne 2$ respectively.
\end{thm}
\demo
Let $G$ be generated by $\{x_1,x_2,\dots,x_{n},y_1,y_2,\dots,y_{n}\}$ with defining relationships $x_i^p=y_i^p=1$ and $[x_i,y_j]=1$ for all $i,j$ such that $1\le i,j\le n$, $[x_{(u-1)a+i},x_{(u-1)a+j}]=1$ for $1\le u \le r$ and $1\le i<j\le a$, $[x_{(u-1)a+i},x_{(v-1)a+j}]=z_{(u-1)a+i,(v-1)a+j}$ for $1\le u<v \le r$ and $1\le i,j\le a$, $[y_i,y_j]=[x_j,x_i]$ for every $i,j$ with $1\le i<j\le n$,
$z_{(u-1)a+1,(v-1)a+j}\in Z(G)$ for $1\le u<v \le r$ and $1\le j\le a$.
 For convenience, we use addition operation to replace mutiplication operation of $G$. We also use the following notations (where $1\le u,v\le r$ and $1\le i,j,k\le a$):
$$Z_{uv}^{(k)}:=(z_{(u-1)a+k,(v-1)a+1},z_{(u-1)a+k,(v-1)a+2},\dots,z_{(u-1)a+k,(v-1)a+a}),$$
$$Z_{uv}=(z_{(u-1)a+i,(v-1)a+j})=\left(\begin{array}{c}
                                       Z_{uv}^{(1)} \\
                                      Z_{uv}^{(2)}\\
                                       \vdots\\
                                      Z_{uv}^{(3)}  \\
                                     \end{array}
                                   \right)\ \mbox{and}\ Z=\left(
                                         \begin{array}{cccc}
                                           Z_{11} & Z_{12} & \dots & Z_{1r} \\
                                            Z_{21} & Z_{22} & \dots & Z_{2r} \\
                                           \vdots &\vdots&\ddots&\vdots \\
                                            Z_{r1} & Z_{r2} & \dots & Z_{rr}\\
                                         \end{array}
                                       \right).$$
Using above notations, we continue to give defining relationships $Z_{uv}^{(k)}=Z_{uv}^{(1)}B^{k-1}$ for $1\le u< v\le r$ and $2\le k\le a$. It is easy to see that $Z_{uv}=-Z_{vu}^T$ for $1\le v<u\le r$,
 $G'=Z(G)=\lg z_{(u-1)a+1,(v-1)a+j}\mid 1\le u<v \le r, 1\le j\le a\rg$ is elementary abelian of order $p^{\frac{r-1}{2}n}$. Hence $|G|=p^{\frac{r+3}{2}n}$.

Since $[x_{(u-1)a+i},x_{(u-1)a+j}]=1$ for $1\le u \le r$ and $1\le i<j\le a$, $Z_{uu}=O_a$ for $1\le u\le r$. Since $Z_{uv}^{(k)}=Z_{uv}^{(1)}B^{k-1}$ for $1\le u< v\le r$ and $2\le k\le a$, $B^TZ_{uv}=Z_{uv}B$ for $1\le u< v\le r$. By Lemma \ref{lem=ss}, $Z_{uv}$ is symmetric for $1\le u< v\le r$. Hence $Z_{uv}=-Z_{vu}^T=-Z_{vu}$ for $1\le v<u\le r$. Moreover $B^TZ_{uv}=Z_{uv}B$ for all $1\le u,v\le r$ and $Z^T=-Z$.
 Let $X=\lg x_1,x_2,\dots,x_{n}\rg Z(G)$ and $Y=\lg y_1,y_2,\dots,y_{n}\rg Z(G)$.

Similar to Assertion 1 in the proof of Theorem \ref{thm=ex1}, we have $|[x,X]|\ge p^{(r-1)a}$ for all $x\in X\setminus Z(G)$. It follows that $|C_X(x)|\le |X|/p^{(r-1)a}=p^a|Z(G)|$.

Similarly, $|C_Y(y)|\le p^{a}|Z(G)|$ for all $y\in Y\setminus Z(G)$.
It is easy to check that $C_G(X)=Y$ and $C_G(Y)=X$, yielding $m_G(G)=m_G(X)=m_G(Y)=p^{(r+1)n}$.

\medskip

Similar to Assertion (1) and Assertion (2) in the proof of Theorem \ref{thm=ex2}, we have $m^*(G)=p^{(r+1)n}$ and $\mathcal{CD}(G)$ is a quasi-antichain.

Let $M$ be an atom different from $X$ and $Y$. Then, by the same reason given in Theorem \ref{thm=main}, we may let $M=\lg w_1,w_2,\dots,w_{n}\rg$ and $N=C_G(M)=\lg v_1,v_2,\dots,v_{3n}\rg$ where
$$(w_1,w_2,\dots,w_{n})=(x_1,x_2,\dots,x_{n})+(y_1,y_2,\dots,y_{n})C,$$
$$(v_1,v_2,\dots,v_{n})=(x_1,x_2,\dots,x_{n})D+(y_1,y_2,\dots,y_{n}).$$
Since $[M,N]=0$, $C^TZ=ZD$. Let
$$C=\left(
                                         \begin{array}{cccc}
                                           C_{11} & C_{12} & \dots & C_{1r} \\
                                            C_{21} & C_{22} & \dots & C_{2r} \\
                                           \vdots &\vdots&\ddots&\vdots \\
                                            C_{r1} & C_{r2} & \dots & C_{rr}\\
                                         \end{array}
                                       \right)\ \mbox{and}\ D=\left(
                                         \begin{array}{cccc}
                                           D_{11} & D_{12} & \dots & D_{1r} \\
                                            D_{21} & D_{22} & \dots & D_{2r} \\
                                           \vdots &\vdots&\ddots&\vdots \\
                                            D_{r1} & D_{r2} & \dots & D_{rr}\\
                                         \end{array}
                                       \right).$$
Similar to Assertion 3 in the proof of Theorem \ref{thm=ex1}, we have $C=D=\Diag(B^k,B^k,\dots,B^k)$ for some $1\le k\le p^a-1$. Hence $w=p^a+1$. It is easy to see that $M=N$ if and only if $CD=I_n$ if and only if $(p^a-1)\mid 2k$. Hence the number of abelian atoms is $1$ or $2$ for $p=2$ and $p\ne 2$ respectively.
\qed

\bigskip
\begin{thm}
  Suppose that $p$ is and odd prime, $a$ is odd, and $r$ is positive integers such that $r\ge 3$. Let $n=ar$. Then there exists a group $G$ such that $|G/Z(G)|=p^{2n}$, and $\mathcal{CD}(G)$ is a quasi-antichain of width $w=p^a+1$, in which the number of abelian atoms is $t=0$.
\end{thm}
\demo Let $G$ be generated by $\{x_1,x_2,\dots,x_{n},y_1,y_2,\dots,y_{n}\}$ with defining relationships similar to that in Theorem \ref{thm=ex3}. The unique different defining relationship is $[y_i,y_j]=[x_j,x_i]^\nu$ for every $i,j$ with $1\le i<j\le n$, where $\nu$ is a fixed quadratic non-residue module $p$.

Similar to the proof of Theorem \ref{thm=ex3}, we have $m^*(G)=p^{(r+1)n}$, $\mathcal{CD}(G)$ is a quasi-antichain, $w=p^a+1$,
$\nu C=D=\Diag(B^k,B^k,\dots,B^k)$ for some $1\le k\le p^a-1$, and $M=N$ if and only if $\nu B^{2k}=I_a$.

By the definition of $B$, $\lg B^{\frac{p^a-1}{p-1}}\rg=\lg lI_a\mid 1\le l\le p-1\rg$. Since $\nu$ is not a square, there exists an odd $m$ such that $B^{m\frac{p^a-1}{p-1}}=\nu I_a$.

If $\nu B^{2k}=I_a$, then $(p^a-1)\mid 2k+m\frac{p^a-1}{p-1}$. Notice that $m$ and $\frac{p^a-1}{p-1}=1+p+\dots+p^{a-1}$ are all odd. There is no integer $k$ such that $\nu B^{2k}=I_a$.
Hence the number of abelian atoms is $t=0$.
\qed


\begin{thebibliography}{99}


\bibitem{ABQW}
 AN, L., Brenna, J., Qu, H., Wilcox, E., Chermak-Delgado lattice extension theorems, {\it Algebra comm.}, 2015, 43(5): 2201--2213.

 \bibitem{BHW1}
 Brewster, B., Hauck, P., Wilcox, E.,  Groups whose Chermak-Delgado lattice is a chain, {\it Journal of Group Theory}, 2014, 17(2): 253--265.

\bibitem{BHW2}
 Brewster, B., Hauck, P., Wilcox, E., Quasi-antichain Chermak-Delgado lattice of finite groups, {\it Arch. Math.}, 2014, 103: 301--311.


 \bibitem{BW}
 Brewster, B., Wilcox, E., Some groups with computable Chermak-Delgado lattices, {\it Bulletin of the Australian Mathematical Society}, 2012, 86(1): 29--40.

 \bibitem{CD}
 Chermak, A.,  Delgado, A.,  A measuring argument for finite groups, {\it Proceedings of the American Mathematical Society}, 1989, 107(4): 907--914.

\bibitem{I} Isaacs, I. M.,  Finite Group Theory,  American Mathematical Society, 2008.

\end{thebibliography}
\end{document}